\documentclass[11pt]{article} \usepackage[dvips]{epsfig}
\usepackage{latexsym}
\usepackage{amsfonts,amsmath,amssymb}
\newtheorem{corollary}{Corollary}

\newtheorem{theorem}{Theorem}
\newtheorem{proposition}{Proposition}

\newcommand{\qed}{\mbox{$\Diamond$}\vspace{\baselineskip}}
\newenvironment{proof}{\noindent{\bf Proof:}}{\qed}

\begin{document}

\author{Mikl\'os B\'ona\\
        Department of Mathematics\\
University of Florida\\
Gainesville FL 32611-8105\\
USA  \and
Istv\'an Mez\H{o}\thanks{The research of Istv\'an Mez\H{o} was supported by the Scientific Research Foundation of Nanjing University of Information Science \& Technology, and The Startup Foundation for Introducing Talent of NUIST. Project no.: S8113062001}\\
        Department of Mathematics\\
Nanjing University of Information Science and Technology\\
Nanjing, 210044, P. R. China}

\title{Real Zeros and Partitions without singleton blocks}
\maketitle

\begin{abstract}
We prove that the generating polynomials of partitions of an $n$-element
set into non-singleton blocks, counted by the number of blocks, have
real roots only and we study the asymptotic behavior of the leftmost roots. We apply this information to find the most likely number
of blocks. 
\end{abstract}

\section{Introduction}
A {\em partition} of the set $[n]=\{1,2,\cdots ,n\}$ is a set of disjoint blocks $B_1,B_2,\cdots ,B_k$ so that $\cup_{i=1}^k B_i=[n]$. 
The number of partitions of $[n]$ into $k$ blocks is denoted by $S(n,k)$ and
is called a {\em Stirling number of the second kind}. 

Similarly, the number of permutations of length $n$ with exactly $k$
cycles is denoted by $c(n,k)$, and is called a {\em signless Stirling number
of the first kind}. See any textbook on Introductory Combinatorics, such
as \cite{walk} or \cite{intro} for the relevant definitions, or basic facts, 
on Stirling numbers. 

The ``horizontal'' generating functions, or generating polynomials, of
 Stirling numbers have many interesting properties. 
Let $n$ be a fixed positive integer. Then it is well-known (see
\cite{walk} or \cite{intro} for instance) that 
\begin{equation} \label{1stir} 
C_n(x)=\sum_{k=1}^n c(n,k)x^k = x(x+1)\cdots (x+n-1).\end{equation}
In particular, the roots of the generating polynomial $C_n(x)$ are all 
real (indeed, they are the integers $0, -1, -2,\cdots ,-(n-1)$). 

Similarly, it is known (see \cite{wilfg}, page 20, for instance) 
that for any fixed positive integer $n$, the
roots of the generating polynomial
\[S_n(x)= \sum_{k=1}^n S(n,k)x^k \] are all real, though they are 
not nearly as easy to describe as those of $C_n(x)$.

Rodney Canfield \cite{canfield} (in the case of $r=1$) and
Francesco Brenti \cite{conbrenti}
(in the general case) have generalized (\ref{1stir}) as follows. Let
$d_r(n,k)$ be the number of permutations of length $n$ that have $k$
cycles, each longer than $r$. Such permutations are sometimes called
$r$-derangements. Then the generating polynomial 
\begin{equation} \label{g1stir}
d_{n,r}(x)=\sum_{k\geq 1} d_r(n,k)x^k 
\end{equation} 
has real roots only. The first author \cite{derangements}
 proved that for any given positive integer constant $m$, there exists
a positive number $N$ so that  if $n>N$, then
one of these roots will be very close to $-1$, one will be very close to
 $-2$, and so on, with one being very close to $-m$, to close
the sequence of $m$ roots being very close to consecutive negative integers. 

In this paper, we consider the analogue problem for set partitions. 
Let $D(n,k)$ be the number of partitions of $[n]$ into $k$ blocks, each
consisting of more than one element. We are going to prove
that the generating polynomial

 \begin{equation} \label{g2stir}
D_n(x)=\sum_{k\geq 1}D(n,k)x^k 
\end{equation} 
has real roots only. We will then use this information 
 to determine the location of the largest coefficient(s) of 
$D_n(x)$. We also prove that the number of blocks
is normally distributed. Finally, we use our methods on $r$-derangements, 
and prove the more general result that for any fixed $r$, the distribution
of the number of cycles
of $r$-derangements of length $n$ converges to a normal distribution.

Note that the fact that the two kinds of Stirling numbers
behave in the same way under this generalization is not completely expected.
Indeed, while $1/e$ of all permutations of length $n$ have no cycles of length
1, and in general, a constant factor of permutations of length $n$ have
no cycles of length $r$ or less, the corresponding statement  is not true
for set partitions. Indeed, almost all partitions of $[n]$ contain a
singleton block as we show in Section 3.1.
 However, as this paper proves, the real zeros property 
survives. 

Finally, we mention that the {\em vertical} generating functions (minimal
block or cycle size is 
fixed, $n$ varies) of permutations and set partitions have been studied in
\cite{canfield}. 

\section{The Proof of The Real Zeros Property}
We start by a recurrence relation satisfied by the numbers
 $D(n,k)$  of partitions of $[n]$ into $k$ blocks, each
block consisting of more than one element. It is straightforward to
see that
\begin{equation}
\label{recu} D(n,k)=kD(n-1,k)+(n-1)D(n-2,k-1).\end{equation}
Indeed, the first term of the right-hand side counts  partitions of 
$[n]$ into blocks larger than one in which the element $n$ is in a block 
larger than two, and the second term of the  right-hand side counts
those in which $n$ is in a block of size exactly two. We note that this recurrence appears in the classical book of Comtet \cite[p. 222]{Comtet}.

Let $D_{n}(x)=\sum_{k\geq 1}D(n,k)x^k$. Then (\ref{recu}) yields
\begin{equation} \label{functioneq}
D_{n}(x)=x\left(D_{n-1}'(x)+(n-1)D_{n-2}(x)\right).
\end{equation}

Note that $D_{1}(x)=0$, and $D_{n}(x)=x$ if $2\leq  n<4$.
So the first non-trivial polynomial $D_{n}(x)$ occurs when $n=4$, and
then $D_4(x)=3x^2+x$. 
In the next non-trivial case of $n=5$, we get $D_5=10x^2
+x$.

\begin{theorem} \label{main} Let $n\geq 2$. Then
the polynomial $D_n(x)$ of degree $\lfloor n/2 \rfloor$ has real roots only. All these roots are simple and non-positive.

Furthermore, the roots of $D_n(x)$ and $D_{n-1}(x)$ are interlacing in the following sense.
If $D_n(x)$ and $D_{n-1}(x)$ are both of degree $d$, and their roots
are, respectively, $0=x_0>x_1>\cdots >x_{d-1}$, and $0=y_0>y_1>\cdots >
y_{d-1}$, then
\begin{equation} \label{samedegree}
0>x_1>y_1>x_2>y_2>\cdots >x_{d-1}>y_{d-1},\end{equation}
while if $D_n(x)$ is of degree $d+1$ and $D_{n-1}(x)$ is of degree $d$,
and their roots are, respectively, $0=x_0>x_1>\cdots >x_{d}$, and
$0=y_0>y_1>\cdots >
y_{d-1}$, then
\begin{equation} \label{diffdegree}
0>x_1>y_1>x_2>y_2>\cdots >x_{d-1}>y_{d-1}>x_{d}.\end{equation}
\end{theorem}

\begin{proof}
We prove our statements by induction on $n$. For $n\leq 4$, the statements are true. Now assume that the statement is true for $n-1$, and let us prove it for $n$. Let $0=y_0>y_1>\cdots >y_{d-1}$ be the roots of $D_{n-1}(x)$.

First we claim that if $0>x>y_1$, then $D_{n-1}(x)<0$, that is,
the polynomial $D_{n-1}$ is negative between its two largest roots.
Indeed, $D_{n-1}'(0)=D(n-1,1)=1$, so $D_{n-1}(x)'>0$
in a neighborhood of 0. This implies that in that neighborhood, $D_{n-1}(x)$ is monotone increasing. As $D_{n-1}(0)=0$, this implies our claim.

Now consider (\ref{functioneq}) at $x=y_1$. We claim that at that root,
we have both $D_{n-1}'(y_1)<0$ and $D_{n-2}(y_1)<0$. The latter is a direct consequence of the previous paragraph and the induction hypothesis. The former follows from the fact that $D_{n-1}(x)<0$ for $y_1<x<0$, the fact that $D_{n-1}(y_1)=0$, and the fact that
the roots of $D_{n-1}$ are all simple by induction.

So when $x=y_1$, the argument of the previous paragraph shows that
the right-hand side of (\ref{functioneq}) is the product
of the negative real number $y_1$, and the negative real number
$D_{n-1}'(y_1)+(n-1)D_{n-2}(y_1)$. Therefore, the left-hand side
must be positive, that is, $D_n(y_1)>0$. As $D_n(x)<0$ in a
neighborhood of 0, this shows that $D_n$ has a root in the
interval $(y_1,0)$.

More generally, we claim that $D_n(x)$ has a root in the interval
$(y_{i+1},y_i)$. For this, it suffices to show that $D_n(y_{i})$ and
$D_n(y_{i+1})$ have opposite signs. This will follow
by (\ref{functioneq}) if we can prove the following two statements.
\begin{enumerate}
\item[(i)] $D_{n-1}'(y_{i})$ and
$D_{n-1}'(y_{i+1})$ have opposite signs ($D_{n-1}'(y_i)$ is negative if and only if $i$ is odd), and

\item[(ii)]
$D_{n-2}(y_{i})$ and $D_{n-2}(y_{i+1})$ have opposite signs, ($D_{n-2}(y_i)$ is negative if and only if $i$ is odd).
\end{enumerate}

Just as before, (i) follows from Rolle's theorem, and (ii) follows from
the induction hypothesis.

As we know that both $D_{n-1}'(y_1)$ and $D_{n-2}(y_1)$ are negative, it is a direct consequence of the preceding two statements that $D_{n-1}'(y_{i})$ and $D_{n-2}(y_{i})$ have equal signs for all $i$.

Therefore, by (\ref{functioneq}), $D_n(y_{i})$ and $D_n(y_{i+1})$ have opposite signs, and so $D_n(x)$ has a root in $(y_{i+1},y_i)$.

The above argument completes the proof of the theorem for odd $n$.

When $n$ is even, then $D_n$ is of degree $d+1$, while $D_{n-1}$ is of degree $d$. In that case, we still have to show that $D_n$ has a root in the interval $(-\infty, y_{d-1})$. However, this follows from the previous paragraph since the last root $x_d$ of $D_n$ must
be negative, and cannot be in any of the intervals $(y_{i+1},y_{i})$.
\end{proof}​

It follows from Theorem \ref{main} that both the sequence
$D_4,D_6,D_8,\cdots $, and the sequence $D_5,D_7,D_9, \cdots $ 
 are {\em Sturm sequences}. The interested reader may consult 
\cite{wilfs} for the definition and properties of Sturm sequences. 

\section{Applications of The Real Zeros Property}

In this Section, we consider two applications of the real zeros property.
Both are combinatorial with a probabilistic flavor. 

\subsection{Locating peaks}

If a polynomial $\sum_{k=1}^n b_kx^k$ with positive coefficients has real
roots only, then it is known \cite{intro} that the sequence $b_1,b_2,\cdots
b_n$ of its coefficients is {\em strongly log-concave}. That is, for
all indices $2\leq j \leq n-1$, the inequality
\[b_j^2 \geq b_{j-1}b_{j+1} \frac{j+1}{j} \cdot \frac{n-j+1}{n-j}\]
holds. In other words, the ratio $b_{j+1}/b_j$ is strictly decreasing
with $j$, and therefore there is at most one index $j$ so that
 $b_{j+1}/b_j=1$. Thus the sequence  $b_1,b_2,\cdots,
b_n$ has either one peak, or two consecutive peaks. 

A useful tool in finding the location of 
this peak is the following theorem of Darroch.

\begin{theorem} \label{darroch} \cite{darroch}
Let $A(x)=\sum_{k=1}^n a_kx^k$ be a polynomial that has
real roots only that satisfies $A(1)>0$. Let $m$ be the index for a peak for the sequence of the coefficients of $A(x)$. Let $\mu=A'(1)/A(1)=
\frac{\sum_{k=1}^nka_k}{\sum_{k=1}^na_k}$.
Then we have
\[|\mu -m| <1 .\] 
\end{theorem}

Note that in a combinatorial setup, $\mu$ is the average value
 of the statistic
 counted by the generating polynomial $A(x)$. For instance, if $A(x)=
S_n(x)$, then $\mu$ is the average number of blocks in a randomly
selected partition of $[n]$. 

There is a very extensive list of results on the peak (or two peaks)
of the sequence $S(n,1), S(n,2), \cdots , S(n,n)$ of Stirling numbers of 
the second kind. See \cite{pom} for a brief history of this topic and the
relevant references.  In particular,  if $K(n)$ denotes the 
index of this peak
(or the one that comes first, if there are two of them), then  $K(n) \sim
n/\log n$ . More precisely, let $r$ be the unique positive root of the
equation 
\begin{equation} \label{defofr}
re^r=n.\end{equation}
 Then, for $n$ sufficiently large, $K(n)$ is one of the
two integers that are closest to $e^r-1$. In view of Theorem \ref{darroch},
 one way to  approach this
problem is by computing the average number of blocks in a randomly selected
partition of $[n]$.  

Now that we have proved that  the generating polynomial 
$D_n(x)=\sum_{k\geq 1} D(n,k)x^k $
has real roots only, it is natural to ask how much of the long list of 
results on Stirling numbers can be generalized to the numbers $D(n,k)$. 
In this paper, we will show a quick way of estimating the average number of
blocks in a partition of $[n]$ with no singleton blocks, and so, 
by Darroch's theorem, the location of the peak(s) in the sequence
$D(n,1), D(n,2), \cdots, D(n,\lfloor n/2 \rfloor)$.
For shortness, let us introduce the notation
$D(n)=\sum_{k}D(n,k)$.

\begin{proposition} \label{pave} Let $X_n$ be the random variable counting
blocks of partitions of $[n]$ that have no singleton blocks. Then
for all positive integers $n\geq 2$, the equality
\begin{equation} \label{average}
E(X_n)=\frac{D(n+1)-n(D(n-1))}{D(n)}
\end{equation}
holds, where $E(X_n)$ denotes the expectation of $X_n$. 
\end{proposition}

\begin{proof} The total number of blocks in all partitions counted by
$D(n)$ is clearly $\sum_{k\geq 1} kD(n,k)$. On the other hand, 
\[\sum_{k\geq 1} kD(n,k)=D(n+1)-nD(n-1),\]
as it follows directly from \eqref{recu}.
\end{proof}

So the peak of the sequence $D(n,1), D(n,2), \cdots $ is one of the two
integers bracketing
\begin{equation}
\frac{D(n+1)}{D(n)}-\frac{n(D(n-1))}{D(n)}.\label{peakbracket}
\end{equation}
We can compare this number with the location  $K(n)$ of the peak of the
 sequence
$S(n,1),S(n,2),\cdots ,S(n,n)$ as follows. 

Let $B(n)$ denote the number of all partitions of $[n]$. This number is
often called a {\em Bell number}. There are numerous precise results on
the asymptotics of the Bell numbers. We will only need the following
fact \cite{debruin}.
\begin{equation} \label{rough} \log B(n) = 
n\left(\log n -\log \log n +O(1)\right),\end{equation}
and its consequence that
\begin{equation} \label{conseq} 
\frac{B(n)}{B(n-1)} \sim \frac{n}{e\log n}.
\end{equation}


For a partition $\pi$ of $[n]$, let $Y_n(\pi)$ be the number of blocks of $\pi$, and let $S_n(\pi)$ be the 
number of singleton blocks of $\pi$.

As the average number of blocks in unrestricted partitions
of $[n]$ is $\frac{1}{B(n)}\sum_{k=1}^n kS(n,k)= \frac{B(n+1)-B(n)}{B(n)}$,
we have  
\begin{equation} \label{stirlingpeak}
E(Y_n)= \frac{B(n+1)}{B(n)}-1 \sim \frac{n}{e\log n}.
\end{equation}

Before comparing formulae (\ref{average}) and (\ref{stirlingpeak}), 
we mention some  simple facts.

For any given element $i\in [n]$, the probability that in a randomly
selected unrestricted partition of $[n]$, the element $i$ forms 
a singleton block is $\frac{B(n-1)}{B(n)}$. Therefore, by linearity
of expectation, we have
\begin{equation} \label{singleton}
E(S_n) = n \frac{B(n-1)}{B(n)} \sim e\log n.
\end{equation}

The following simple result will be very useful for us, and
therefore, we state it as a proposition.

\begin{proposition} \label{simpleprop}
For all positive integers $n$, the equality
\[B(n)=D(n)+D(n+1) \] holds.
\end{proposition}

\begin{proof}  We define a simple  bijection $f$ from the set of partitions
of $[n]$ and $[n+1]$ with no singleton blocks into the set of partitions
of $[n]$. 
On partitions counted by $D(n)$, let $f$ act as the identity map.
On partitions counted by $D(n+1)$, let $f$ act by removing
the element $n+1$ and turning each element that shared a block with $n+1$ into
a singleton block.

\end{proof}

This simple fact has two important corollaries that we will use.
\begin{corollary}\label{DasympB}
We have  $D(n+1)\sim B(n)$.
\end{corollary}

\begin{proof}
Note that $\frac{B(n)}{D(n)} \rightarrow \infty$ since $D(n+1)<B(n)$ and
$\frac{D(n+1)}{D(n)} \rightarrow \infty $. To see the latter, note that
for any $h$, there exists an $N$ so that if $n>N$, then almost all 
partitions counted by $D(n)$ have more than $h$ blocks.

As adding the entry $n+1$  to any block of any partition counted by $D(n)$ results in a partition
counted by $D(n+1)$, the inequality $\lim_{n\rightarrow \infty} D(n+1) / D(n) \geq h $ follows for any $h$.
\end{proof}

Now we can easily see that the locations of the peaks of the sequences $D(n+1,1), 
D(n+1,2),
\cdots $, and $B(n,1), B(n,2), \cdots $, as well as the averages
$E(X_{n+1})$ and $E(Y_n)$ as given in 
 formulae (\ref{average}) and (\ref{stirlingpeak}) are indeed very close to
each other. 
In fact, $E(X_{n+1})\sim E(Y_n)$ as can be seen by comparing
 (\ref{average}) and (\ref{stirlingpeak}). We will not attempt a more precise
comparison here. However, we would like to point out that 
the $-\frac{(n+1)(D(n)}{D(n+1)}$ summand in (\ref{average}), when
$n$ is replaced by $n+1$, asymptotically agrees with $E(S_n)$ as computed
in (\ref{singleton}).  This is in line with what one would intuitively expect,
since the difference between partitions on which $X_{n+1}$ is defined
and partitions on which $Y_n$ is defined is that in the former, singleton
blocks are not allowed.

\begin{corollary} \label{relation}
Let $N_n$ be the variable counting the non-singleton blocks of a randomly
selected unrestricted partition of $[n]$. Let $X_n$ denote the number of
blocks of a randomly selected partition of $[n]$ with no singleton blocks.

Then we have
\begin{equation} \label{firstcor}
E(X_{n+1}-1)\frac{D(n+1)}{B(n)} + E(X_n)\frac{D(n)}{B(n)}  =  E(N_n),
\end{equation}
and also,
\begin{equation}
E((X_{n+1}-1)^2) \frac{D(n+1)}{B(n)} + E(X_n^2)\frac{D(n)}{B(n)} =E(N_n^2).
\end{equation}
\end{corollary}

\begin{proof} Direct consequence of the bijection
$f$ defined in the proof of Proposition \ref{simpleprop}. 
\end{proof}

\subsection{Another way to locate the peaks}

From the results of the previous section -- see \eqref{peakbracket}, \eqref{rough}, \eqref{conseq} and Corollary \ref{DasympB} -- it follows that the asymptotic location of the $K_n^*$ peak of the sequence $D(n,1), D(n,2),\dots,D(n,n)$ is
\begin{equation}
K_n^*\sim\frac{n}{\log(n)},\label{peakasymp}
\end{equation}
which is the same as for the classical Stirling numbers \cite{harper,menon,rendob}.

Analyzing the ordinary generating function of $D(n,k)$ with the saddle point method we are going to show that as $n$ goes to infinity,
\begin{equation}
D(n,k)\sim S(n,k)\sim\frac{k^n}{k!}\label{Dnkasymp}
\end{equation}
for any fixed $k$. From this it will follow at once that the asymptotics \eqref{peakasymp} holds. What is more, relying on \eqref{Dnkasymp} and following the proof presented in \cite{mc} and in its references one can prove that the maximizing index is close to
\[\frac{n-\frac12}{W\left(n-\frac12\right)},\]
where $W(n)$ is the Lambert function and it is the unique solution of the equation $W(n)e^{W(n)}=n$.

To prove \eqref{Dnkasymp} we need the following proposition.

\begin{proposition} \label{ordgf}
For any fixed integers $k$, let
\[f_k(x)=\sum_{n=0}^\infty D(n,k)x^n\]
be the ordinary generating function of the sequence $(D(n,k))$. The following recursion holds true:
\[f_k(x)=\frac{x^2}{1-kx}\left(xf_{k-1}(x)\right)', \quad(k\ge2)\]
with
\[f_1(x)=\frac{x^2}{1-x}.\]
\end{proposition}

\begin{proof}
From \eqref{recu} it follows that
\[f_k(x)=\sum_{n=0}^\infty D(n,k)x^n=\]
\[kx\sum_{n=0}^\infty D(n-1,k)x^{n-1}+x^2\sum_{n=2}^\infty (n-1)D(n-2,k-1)x^{n-2}=\]
\[kx\sum_{n=0}^\infty D(n,k)x^n+x^2\sum_{n=2}^\infty (n+1)D(n,k-1)x^n,\]
which is equivalent to our recursion. The form of $f_1(x)$ is obvious, since $D(n,1)=1$ if $n\ge2$.
\end{proof}

This simple observation and induction shows that, in general, the $f_k(x)$ functions are rational functions of the form
\begin{equation}
f_k(x)=x^{2k}\frac{p_k(x)}{(kx-1)((k-1)x-1)^2((k-2)x-1)^3\cdots(x-1)^k}\quad(k\ge1),\label{genformoff}
\end{equation}
where $p_k(x)$ is a polynomial of degree $\frac{k(k-1)}{2}$. These polynomials first appeared in a 1934 paper of Ward \cite{ward} who studied the representations of the classical Stirling numbers as sums of factorials (see \cite{clark} for more details and other citations). The first $f_k(x)$ functions are as follows:
\begin{align*}
f_1(x)&=x^2\frac{-1}{x-1}\\
f_2(x)&=x^4\frac{2x-3}{(x-1)^2(2x-1)}\\
f_3(x)&=x^6\frac{-12 x^3+40 x^2-45 x+15}{(x-1)^3 (2x-1)^2 (3x-1)}\\
f_4(x)&=x^8\frac{ 288 x^6-1560 x^5+3500 x^4-4130 x^3+2625 x^2-840 x+105}{(x-1)^4 (2 x-1)^3 (3 x-1)^2 (4 x-1)}.
\end{align*}

Going back to our original goal, formula \eqref{genformoff} enables us to prove the following.
\begin{proposition}\label{asympt}For any fixed positive integer $k$ and large $n$ we have that
\[D(n,k)=\frac{k^n}{k!}+O(k-1+\varepsilon)^n\]
holds for arbitrary $\varepsilon>0$.
\end{proposition}

\begin{proof}
The asymptotics of $D(n,k)$ can be determined by analyzing the singularities of its generating function. This is the well known saddle point method described in details by Wilf in \cite{wilfg}.

The function $f_k(x)$ has $k$ singular points on the real line and the smallest one is at $x_0=\frac{1}{k}$. This pole is of order one. The principal part of $f_k(x)$ around this point is
\[PP\left(f_k,\frac{1}{k}\right)=-\frac{1}{k\cdot k!\left(x-\frac{1}{k}\right)}.\]
Since $x_0$ is a first order pole and there are no more poles with the same absolute value, the saddle point method \cite[Theorem 5.2.1]{wilfg} in this particular case says that
\[D(n,k)=[x^n]PP\left(f_k,\frac{1}{k}\right)+O\left(\frac{1}{R'}+\varepsilon\right)^n,\]
where $R'$ is the modulus of the second smallest singular point. In this case this point is $x_1=\frac{1}{k-1}$. Expanding the $PP\left(f_k,\frac{1}{k}\right)$ principal part with respect to $x$ we get the statement.
\end{proof}

\subsection{The asymptotics of the zeros of $D_n(x)$}

Having proven that the zeros of the $D_n(x)$ polynomials are all real (and negative), it can be asked that {\em how large} is the leftmost zero of $D_n(x)$? Let $z_n^*$ denote this {\em leftmost} zero. We point out that an easily calculable upper bound can be given, and this upper bound   approximates $z_n^*$ surprisingly well. This approximation is based on a theorem of Laguerre and Samuelson.

Let
\begin{equation}
p(x)=x^n+a_1x^{n-1}+\cdots+a_{n-1}x+a_n\label{pxgen}
\end{equation}
be an arbitrary polynomial such that $a_n\neq0$. Samuelson \cite{samuelson} -- rediscovering the results of Laguerre \cite{laguerre} -- gave bounds for the interval which contains all the zeros of a polynomial \eqref{pxgen} whose zeros are all \emph{real}. Samuelson's result states that all these zeros are contained in the interval $[x_-,x_+]$, where
\begin{equation}
x_{\pm}=-\frac{a_1}{n}\pm\frac{n-1}{n}\sqrt{a_1^2-\frac{2n}{n-1}a_2}\label{Sam}
\end{equation}
for \eqref{pxgen}.

We want to determine $x_-$ when $p(x)=D_n(x)$ (obviously, $x_+=0$). It can easily be seen that
\[a_1=\frac{D\left(n,\lfloor\frac{n}{2}\rfloor-1\right)}{D\left(n,\lfloor\frac{n}{2}\rfloor\right)},\quad\mbox{and}\quad a_2=\frac{D\left(n,\lfloor\frac{n}{2}\rfloor-2\right)}{D\left(n,\lfloor\frac{n}{2}\rfloor\right)}.\]
Hence
\[|z_n^*|\le\frac{D\left(n,\lfloor\frac{n}{2}\rfloor-1\right)}{\lfloor n/2\rfloor D\left(n,\lfloor\frac{n}{2}\rfloor\right)}+\frac{\lfloor n/2\rfloor-1}{\lfloor n/2\rfloor}\sqrt{\left(\frac{D\left(n,\lfloor\frac{n}{2}\rfloor-1\right)}{D\left(n,\lfloor\frac{n}{2}\rfloor\right)}\right)^2-\frac{2\lfloor n/2\rfloor}{\lfloor n/2\rfloor-1}\frac{D\left(n,\lfloor\frac{n}{2}\rfloor-2\right)}{D\left(n,\lfloor\frac{n}{2}\rfloor\right)}}.\]
In the particular case of the $D_n(x)$ polynomials the Samuelson estimation works very well. The following table compares the actual values of $|z_n^*|$ with the estimates obtained by the Laguerre-Samuelson theorem for 
some even numbers $n$. (Note that $\deg D_n(x)=\lfloor n/2\rfloor$.)
\begin{center}
\begin{tabular}{c|c|c|c}
$n$&10&100&200\\
Numerical value of $|z_n^*|$&9.22&11\,085.5&89\,380.6\\
Estimate of Samuelson&9.24&11\,160.8&90\,011.4
\end{tabular}
\end{center}
The following table contains the analogous information for  odd numbers $n$.
\begin{center}
\begin{tabular}{c|c|c|c}
$n$&11&101&201\\
Numerical value of $|z_n^*|$&2.828&2\,852.96&22\,677.2\\
Estimate of Samuelson&2.85&2\,958.05&23\,552.4
\end{tabular}
\end{center}

By simple combinatorial arguments one can find the special values of $D\left(n,\lfloor\frac{n}{2}\rfloor\right)$, $D\left(n,\lfloor\frac{n}{2}\rfloor-1\right)$, and $D\left(n,\lfloor\frac{n}{2}\rfloor-2\right)$ easily. For example,
\[D\left(n,\left\lfloor\frac{n}{2}\right\rfloor\right)=\left\{\begin{array}{ll}\frac{n!}{2^{n/2}\left(\frac{n}{2}\right)!},&\mbox{if $n$ is even};\\\binom{n}{3}\frac{(n-3)!}{2^{(n-3)/2}\left(\frac{n-3}{2}\right)!},&\mbox{if $n\ge3$ is odd}.\end{array}\right.\]
With these special values one can find that for the Samuelson estimate
\begin{align*}
x_-&\sim-\frac{1}{36\sqrt{6}}n^3\quad(\mbox{$n\to\infty$ is even}),\\
x_-&\sim-\frac{1}{108\sqrt{10}}n^3\quad(\mbox{$n\to\infty$ is odd}).
\end{align*}
These asymptotics and the numerical calculations suggest the following conjecture about the asymptotic behavior of the leftmost zero of $D_n(x)$:
\begin{align*}
z_n^*&\sim-c_{\mathrm{even}}n^3\quad(\mbox{$n\to\infty$ is even}),\\
z_n^*&\sim-c_{\mathrm{odd}}n^3\quad\;(\mbox{$n\to\infty$ is odd}).
\end{align*}

\section{Basic modularity properties of $D(n,k)$ and $D_n(1)$}

A simple application of the binomial theorem (see \eqref{bintrfid} below) reveals that the well known modularity property
\begin{equation}
S(p,k)\equiv0\pmod{p}\quad(1<k<p)\label{stdiv}
\end{equation}
of the Stirling numbers can be transferred to the $D(n,k)$ numbers.

\begin{proposition}\label{divprop}
\begin{equation}
D(p,k)\equiv0\pmod{p}\quad(1<k<p).\label{ddiv}
\end{equation}
Here $p$ is an arbitrary prime.
\end{proposition}
The most basic Bell number divisibility follows directly from \eqref{stdiv}:
\[B(p)\equiv2\pmod{p}\]
for odd primes $p$. The corresponding divisibility for $D_n=D_n(1)$ is the consequence of \eqref{ddiv} and of \eqref{g2stir}:
\begin{corollary}We have that
\[D_p\equiv1\pmod{p}\]
for any prime $p$ including $p=2$.
\end{corollary}

\begin{proof} (of Proposition \ref{divprop}). The identity from which we can deduce the proposition reads as
\begin{equation}
D(n,k)=\sum_{s=n-k}^n\binom{n}{s}(-1)^{n-s}S(s,s+k-n).\label{bintrfid}
\end{equation}
From this and from the fact that $p\left|\binom{p}{k}\right.$ $(1\le k<p)$ \eqref{ddiv} follows, indeed.

Formula \eqref{bintrfid} can be proven easily considering the exponential generating function
\begin{equation}
\sum_{n=0}^\infty D(n,k)\frac{x^n}{n!}=\frac{1}{k!}\left(e^x-1-x\right)^k.\label{egfDnk}
\end{equation}
Expanding $\left((e^x-1)-x\right)^k$ with the binomial theorem and using the fact that $\frac{1}{k!}(e^x-1)^k$ is the exponential generating function of the Stirling numbers, we are done.
\end{proof}

We remark that \eqref{bintrfid} can be generalized to
\[D_m(n,k)=\sum_{s=0}^{n}\binom{n}{s}\sum_{i=0}^k(-1)^{k-i}S(s,i)D_{m-1}(n-s,k-i),\]
where $D_m(n,k)$ counts the partitions of $[n]$ into $k$ blocks such that all of the blocks contain at least $m$ elements (so $D(n,k)=D_2(n,k)$).

Finally, we also note that from \eqref{egfDnk} we get the exponential generating function of the $D_n=D_n(1)$ numbers:
\[\sum_{n=0}^\infty D_n\frac{x^n}{n!}=e^{e^x-1-x},\]
from which a simple pair of formulas follows connecting the Bell numbers and the $D_n$ sequence:
\begin{align*}
D_n&=\sum_{j=0}^n\binom{n}{j}(-1)^jB(n-j),\\
B(n)&=\sum_{j=0}^n\binom{n}{j}D_j.\\
\end{align*}
The combinatorial proofs can be easily given by the inclusion-exclusion principle.

\section{Further Directions}
It is natural to ask whether Theorem \ref{main} can be generalized to 
partitions with all blocks larger than $r$, where $r$ is a given
positive integer. That would parallel the result (\ref{g1stir}) 
of Brenti \cite{conbrenti}  on permutations.

A consequence of the fact that the polynomials $D_n(x)$ have real zeros
is that for any fixed $n$, the sequence $D(n,1), D(n,2),\cdots ,
D(n,\lfloor n/2 \rfloor )$ is log-concave. In \cite{sagan}, Bruce Sagan
provides a proof for the analogous statement for the Stirling numbers of the second kind (i.e., when $r=0$). However, his injection proving
that result does not preserve the no-singleton-block property.
 Now that we know that the statement is
true, it is natural to ask for an injective proof. Similarly, as we know
 that $d_{n,r}(x)$ has real roots only
(see (\ref{g1stir})), 
we can ask for a combinatorial proof for the fact that  
the sequence $d_{r}(n,1), d_r(n,2),\cdots ,
d_r(n,\lfloor n/r \rfloor ) $ is log-concave.

\end{document}